\documentclass[11pt]{amsart}
\usepackage{amsmath}
\usepackage{amsfonts}
\usepackage{amssymb}
\usepackage{latexsym}
\usepackage{epsfig}
\usepackage{graphicx}
\usepackage{oldgerm}
\usepackage{amsthm}
\usepackage{tikz}
\usepackage{hyperref} 
\usepackage{color}

\theoremstyle{theorem}
\newtheorem{thm}{Theorem}[section]
\newtheorem{lem}[thm]{Lemma}
\newtheorem{cor}[thm]{Corollary}
\newtheorem{prop}[thm]{Proposition}

\theoremstyle{remark}
\newtheorem{rem}[thm]{Remark}
\newtheorem{ex}[thm]{Example}

\def\la{\lambda}

\def\zbar{\overline{z}}
\def\wbar{\overline{w}}
\newcommand{\ii}{\mathbf{i}}

\def\res{\mathrm{Res}}
\def\gammabar{\bar{\gamma}}
\def\fpv{f_{\mathrm{PV}}}
\def\C{\mathbb{C}}
\def\D{\mathbb{D}}
\def\E{\mathbb{E}}

\def\Q{\mathbb{Q}}
\def\R{\mathbb{R}}
\def\Z{\mathbb{Z}}
\def\cJ{\mathcal{J}}
\def\cP{\mathcal{P}}
\def\cZ{\mathcal{Z}}

\makeatletter
\def\filename{\texttt{\jobname.tex}} 
\makeatother

\begin{document}
\title[Zeros of random power series with finitely
dependent Gaussian coefficients]{Expected number of zeros of
random power series with finitely dependent Gaussian coefficients}

\author{Kohei Noda} 
\address{Graduate School of Mathematics, Kyushu
University, 744 Motooka, Nishi-ku, Fukuoka 819-0395, Japan} 
\email{noda.kohei.721@s.kyushu-u.ac.jp} 
\author{Tomoyuki Shirai}
\address{Institute of Mathematics for Industry, Kyushu
University, 744 Motooka, Nishi-ku, Fukuoka 819-0395, Japan}
\email{shirai@imi.kyushu-u.ac.jp}
 
 
\begin{abstract}
We are concerned with zeros of random power series with coefficients being a
 stationary, centered, complex Gaussian process. 
We show that the expected number of zeros in every smooth domain 
in the disk of convergence is less than that of the hyperbolic GAF with
 i.i.d. coefficients. 
When coefficients are finitely dependent, i.e., 
the spectral density is a trigonometric polynomial, 
we derive precise asymptotics of the expected number of
 zeros inside the disk of radius $r$ centered at the origin 
as $r$ tends to the radius of convergence, 
in the proof of which we clarify that the negative
 contribution to the number of zeros stems
 from the zeros of the spectral density. 
\end{abstract}
\maketitle

\section{Introduction}

Let $\{\zeta_k\}_{k=0}^{\infty}$ be independent, 
identically distributed (i.i.d.) standard complex Gaussian random variables. 
Peres and Vir\'ag studied the zeros of random power series
$\fpv(z)=\sum_{k=0}^\infty \zeta_k z^k$ and found that the
zero point process $\sum_{z \in \C : \fpv(z)=0}\delta_z$
becomes a determinantal point process associated with the
Bergman kernel \cite{YV}. 
The studies around this Gaussian analytic function (GAF)
has been developing in several directions
(cf. \cite{BH,JBANRPMS,GS, KS, KRI2,SMTS}), however, it seems that there are relatively few works on zeros of
random power series with \textit{dependent} Gaussian
coefficients. 
Recently, Mukeru, Mulaudzi, Nazabanita and Mpanda studied
the zeros of Gaussian random power series $f_H(z)$ on
the unit disk with coefficients $\Xi^{(H)} = \{\xi_k^{(H)}\}_{k=0}^{\infty}$ being a fractional Gaussian noise (fGn) with Hurst index $0\leq
H<1$. They gave an estimate for the expected number of 
zeros of $f_H(z)$ inside $\D(r):=\{z\in\C:|z|<r\}$ and show that it is smaller than that of $\fpv(z)$
by $O((1-r^2)^{-1/2})$ \cite{SMJM}, whose proof was based on the
 maximum principle via an integral representation on $\D(r)$ of the expectation. 
In this paper, we will give a precise asymptotics as $r \to
1_-$ of the expected number of zeros in $\D(r)$ of a random
power series $f_{\Xi}(z) = \sum_{k=0}^{\infty} \xi_k z^k$ when $\Xi = \{\xi_k\}_{k=0}^{\infty}$ is a
stationary, centered, 
finitely dependent complex Gaussian process, i.e., its spectral density is a trigonometric polynomial
of degree $n$. 
As will be seen later, the essential idea of our proof is to represent 
the expected number of zeros as a contour integral on
$\partial \D(r)$ by using the Stokes theorem
similar to \cite{JB,KN} and keep track of  
the poles of the integrand indexed by $r$, i.e., the zeros of a 
(scaled) spectral density for $\Xi$, as $r \to 1_-$. 
We found that the degeneracy of zeros of spectral density
sensitively affects on the order of the difference between 
the expected number of zeros of $f_{\Xi}(z)$ and that of $\fpv(z)$. 

Let $\Xi = \{\xi_k\}_{k \in \Z}$ is a stationary, centered, 
complex Gaussian process with \textit{unit variance} and covariance function 
\begin{equation}
\E[\xi_k\overline{\xi_l}]=\gamma(l-k), \quad k, l \in \Z, 
\label{eq:covGAF}
\end{equation}
where $\gamma(0)=1$ and $\gamma(-k) =
\overline{\gamma(k)}$. 
Throughout this paper, we always assume the variance to be $1$. 
We consider the following random power series 
\begin{equation}
f_{\Xi}(z) = \sum_{k=0}^\infty \xi_k z^k. 
\label{eq:GAF}
\end{equation}
For the sake of simplicity, in what follows, we often omit the subscript $\Xi$
in $f_{\Xi}$. The covariance matrix of the Gaussian analytic function
(GAF) defined in \eqref{eq:GAF} is given by 
\begin{equation}
K_f(z,w) 
= \E[f(z)\overline{f(w)}] = \frac{1}{1-z\wbar} G_2(z,w), 
\label{eq:general_covariance}
\end{equation}
where 
\begin{equation}
G_2(z,w) = 1+G(z) + \overline{G(w)}, \quad 
 G(z) = \sum_{k=1}^\infty \overline{\gamma(k)} z^k. 
\label{eq:Gzw}
\end{equation}
Since $|\gamma(k)| \le \gamma(0)=1$ follows from positive definiteness, 
the convergence radius of $G(z)$ is more than or equal to $1$. The covariance function $\gamma(k)$ can be represented as $\gamma(k) =
(2 \pi)^{-1} \int_0^{2\pi} e^{\sqrt{-1}k\theta} d\Delta(\theta)$, 
where $\Delta(\theta)$ is
called the spectral function of $\Xi$. When $\Delta(\theta)$ is
absolutely continuous with respect to the Lebesgue measure, 
the density $\Delta'(\theta) = d\Delta(\theta)/d\theta$ is called the
spectral density of $\Xi$ (cf. \cite{DM}). We note that
$G_2(e^{\sqrt{-1}\theta}, e^{\sqrt{-1}\theta})$ gives the spectral density of the Gaussian
process $\Xi$ if $G(z)$ is analytic in a neighborhood of
$\D$. 
When $\{\xi_k\}_{k \in \Z}$ are i.i.d., $\gamma(k) =
\delta_{0,k}$ (Kronecker's delta) and $K_f(z,w)$ is the Szeg\H{o} kernel. 
As mentioned before, Peres-Vir\'ag showed that the zeros of $\fpv(z)$
with i.i.d. Gaussian coefficients 
form the determinantal point process associated with the Bergman
kernel \cite{YV}. 
In the present paper, we compare the expected number of
zeros of $f(z)$ with \textit{finitely dependent} Gaussian coefficients with that of
$\fpv(z)$. 

We first deal with the case of $2$-dependent stationary Gaussian
processes with covariance function 
\begin{equation}
  \gamma(k) = 
\begin{cases}
    1 & k=0, \\
    a & |k|= 1, \\
    b & |k|= 2, \\
    0 & \text{otherwise}.  
\label{eq:gammapositive}
  \end{cases}
\end{equation}
We easily verify that $\{\gamma(k)\}_{k \in \Z}$ is positive
definite if and only if $(a,b)$ is in the region 
$\cP = \cP_1 \cup \cP_2$ with 
\begin{equation*}
\cP_1 = \left\{(a,b) \in \R^2 :
	 \frac{a^2}{8}+\left(b-\frac{1}{4}\right)^2 \leq
 \frac{1}{16} \right\}
\end{equation*}
and 
\begin{equation*}
\cP_2 = \left\{ 
(a,b) \in \R^2 : 
\frac{a^2}{8}+\left(b-\frac{1}{4}\right)^2\geq \frac{1}{16},
\ |a|-\frac{1}{2} \le b \le \frac{1}{6}
\right\}. 
\end{equation*}
\begin{figure}[htbp]
\begin{center}
\includegraphics[width=0.4\hsize]{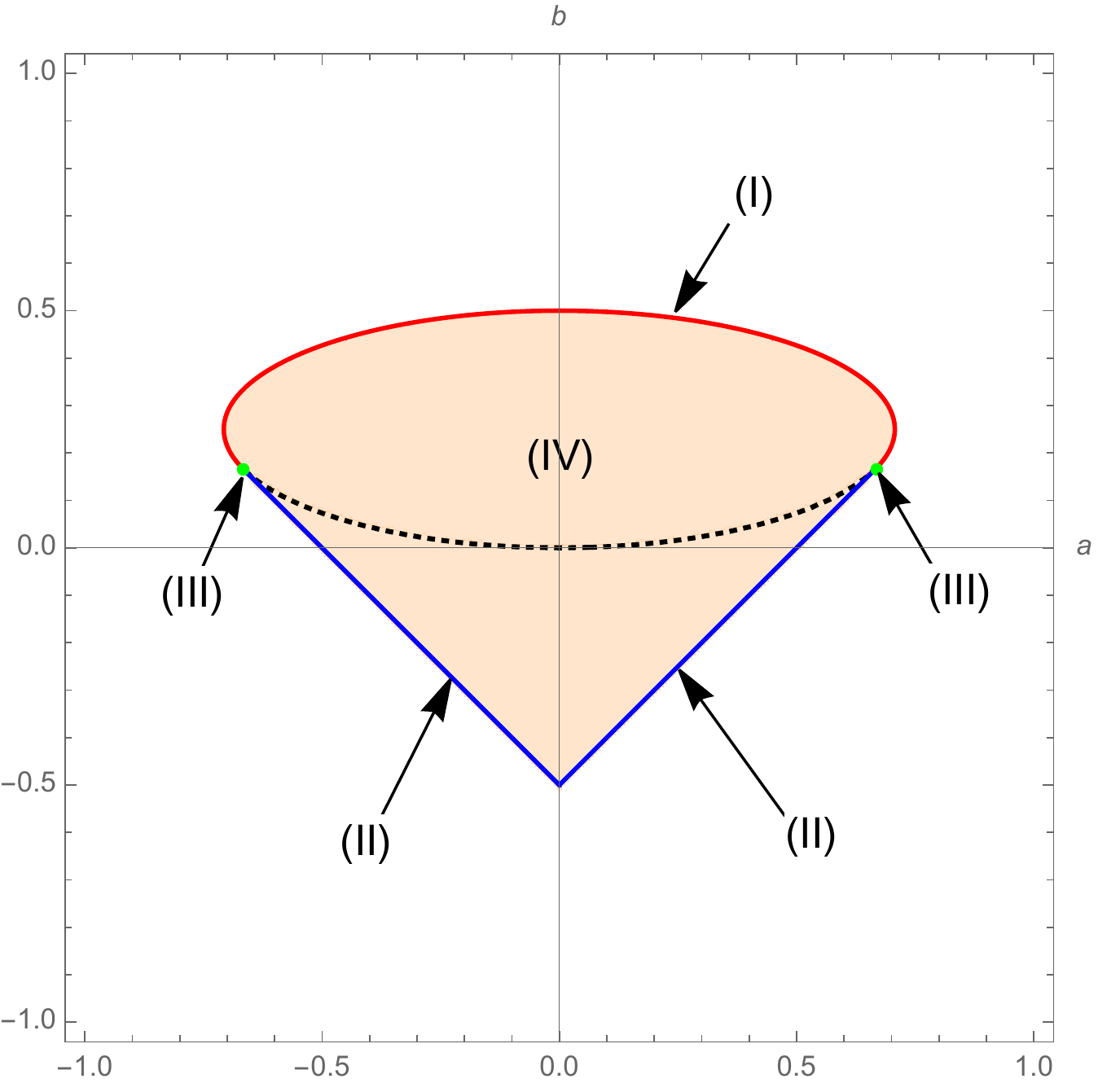}
\end{center}
\caption{The region $\cP$ of positive definiteness of $\gamma(k)$ defined in 
\eqref{eq:gammapositive}. The red and black dashed ellipse is
 the boundary of $\cP_1$, the green points are $(a,b) = (\pm
 2/3, 1/6)$, and the blue line segments is $b=|a|-1/2$ for $-1/2 \le
 b \le 1/6$. 
 } 
\label{fig:N2abintro} 
\end{figure}
See Figure~\ref{fig:N2abintro}. 
We consider the GAF $f_{a,b}(z)$ associated with \eqref{eq:gammapositive}. 
Since we normalized the variance of $\xi_k$ to be $1$, 
the convergence radius of the power series $f_{a,b}(z)$ is
$1$ a.s. for any $(a,b) \in \cP$.  

We denote the zeros of GAF $f$ by $\cZ_f$ and let 
\begin{align*}
 N_f(r)=\#\{z\in \cZ_f : |z|< r\}, \quad r \in (0,1)
\end{align*}
be the number of zeros within $\D(r)$, the disk of radius
$r$ centered at the origin. 
From now on, for simplicity, we write $r \to 1$ instead of $r\to1_-$.

\begin{thm}\label{thm:mainintro}
Let $f_{a,b}$ be the GAF defined in \eqref{eq:GAF} with 
covariance function of the form \eqref{eq:gammapositive}
 with $(a,b) \in \cP$. 
Then the asymptotic behavior of the expected number of zeros is as follows. \\
\upshape{(I)} If $(a,b)$ satisfies $a^2/8+(b-1/4)^2=1/16$ and
 $1/6 < b \leq 1/2$, then  
\begin{equation}
\E N_{f_{a,b}}(r)=\frac{r^2}{1-r^2}-\sqrt{\frac{2b}{6b-1}}\frac{1}{(1-r^2)^{1/2}}+O(1),
 \quad r\to 1.
\label{eq:case1} 
\end{equation}
\upshape{(II)} If $(a,b)$ satisfies $b=|a|-1/2$ and $-1/2\leq
 b<1/6$, then  
\begin{equation}
\E N_{f_{a,b}}(r) =
 \frac{r^2}{1-r^2}-\frac{1}{2}\sqrt{\frac{1-2b}{1-6b}} 
\frac{1}{(1-r^2)^{1/2}}+O(1), \quad r \to 1.
\label{eq:case2}
\end{equation}
(III) If $(a,b)=(\pm2/3,1/6)$, then  
\begin{equation}
\E N_{f_{a,b}}(r) = 
\frac{r^2}{1-r^2}-\frac{1}{2^{5/4}}\frac{1}{(1-r^2)^{3/4}}+O\left(\frac{1}{(1-r^2)^{1/4}}\right),
\quad r\to1.
\label{eq:case3}
\end{equation}
(IV) If $(a,b)$ is in the interior of $\cP$, then there
 exists a non-negative constant $C(a,b)$ such that 
\begin{equation}
\E N_{f_{a,b}}(r) = \frac{r^2}{1-r^2} - C(a,b)+O\left(1-r^2\right),
 \quad r \to 1.
\end{equation}
The constant $C(a,b)$ is positive except for $(a,b)=(0,0)$. 
The numbers (I)--(IV) in Theorem~\ref{thm:mainintro}
correspond to those in Figure~\ref{fig:N2abintro}.  
\end{thm}

The case of $(a,b)=(0,0)$ corresponds to the case of
Peres-Vir\'ag, $\fpv(z)$, and 
it is known that 
\begin{align}
 \E N_{f_{0,0}}(r) = \E N_{\fpv}(r) = \frac{r^2}{1-r^2}.  \label{eq:Nf} 
\end{align}
Therefore, for all cases, the expected number of zeros is less than
that of $\fpv(z)$ at least in the limit as $r \to 1$. 
In fact, we can show the following stronger result. 
\begin{thm}\label{thm:negative}
Let $f$ be a GAF defined in \eqref{eq:GAF} with
 \eqref{eq:general_covariance} and \eqref{eq:Gzw}. 
Let $D \subset \D$ be a domain with smooth boundaries and 
$N_f(D)$ be the number of zeros of $f$ inside $D$. Then, 
$\E N_f(D)$ is always less than or equal to $\E N_{\fpv}(D)$. 
Moreover, the equality holds for some (hence any) domain $D$ if and only if $f(z)$ is equal to
  $\fpv(z)$ in law. 
 \end{thm}

As was seen in the above, the asymptotic behavior at $(a,b)
= (\pm 2/3, 1/6)$ corresponding to Case (III) is special since $G_2(z,z)$ 
is the most degenerated in the sense that 
\[
 G_2(z,z) = 1 \pm \frac{2}{3} (z+z^{-1}) + \frac{1}{6}
 (z^2+z^{-2}) = \frac{1}{6} z^{-2} (z \pm 1)^4
\]
for $z \in \partial\D = \{z \in \C : |z| = 1\}$. 
The above $G_2(z,z)$ has the degenerated zero at $z=\mp 1$. 
The phenomena are the same in both cases and so we only deal
with the $+$ case below. Now we focus on the $n$-dependent
stationary Gaussian process $\Xi$ with covariance function $\{\gamma_n(k)\}_{k \in \Z}$ 
which is the most degenerated in the sense above, i.e., 
\begin{equation}
 \gamma_n(k) = 
\begin{cases}
 {2n \choose n+k} {2n \choose n} ^{-1} & \text{if $|k| =0, 1, 2, \dots, n$}, \\
0 & \text{otherwise}, 
\end{cases}
\label{eq:gamman} 
\end{equation}
which is normalized as $\gamma_n(0)=1$. 
It is easy to see that 
\begin{equation}
 G_2(z,z) 
= \sum_{k = -n}^n \gamma_n(k) z^k = {2n \choose n}^{-1} z^{-n} (z+1)^{2n}
\label{eq:G2degenerated} 
\end{equation}
for $z \in \partial\D$ and $z=-1$ is the zero of order $2n$. 
We remark that for this Gaussian process $\Xi$ we have the following moving-average representation: 
\[
\xi_k = {2n \choose n}^{-1/2} \sum_{j=0}^{n} {n \choose j} \zeta_{k-j}, \quad k=0,1,\dots, 
\]
where $\{\zeta_j\}_{j \in \Z}$ is an i.i.d. standard complex Gaussian sequence. 
In this case, we have the following asymptotics, 
which include \eqref{eq:case3} as a special case of $n=2$. 

\begin{thm}\label{thm:mainintro2}
Let $\gamma_n(k)$ be defined as \eqref{eq:gamman} and 
$\Xi = \{\xi_k \}_{k \in \Z}$ be the
 stationary, centered, complex Gaussian process with covariance
 function $\{\gamma_n(k)\}_{k \in \Z}$. 
The expected number of zeros of the power series $f$ with
 coefficients $\Xi$ within $\D(r)$ is given by 
\begin{equation}
 \E N_f(r) = \frac{r^2}{1-r^2} - D_n
 (1-r^2)^{-\frac{2n-1}{2n}} 
+O((1-r^2)^{-\frac{2n-3}{2n}}), 
\quad r \to 1, 
\label{eq:case3_asymptotics} 
\end{equation}
where 
\[
D_n= \frac{1}{2n \sin\frac{\pi}{2n}} 
\left\{\binom{2(n-1)}{n-1}\right\}^{\frac{1}{2n}}.  
\]
\end{thm}

\begin{rem}
The term of order $(1-r^2)^{-\frac{2n-2}{2n}}$ in \eqref{eq:case3_asymptotics} vanishes 
 by a cancellation. See the proof of Theorem~\ref{thm:mainintro2} and 
Remark~\ref{rem:cancel}. 
\end{rem}

 As will be seen in the proof of the theorems, 
the order of the second term 
in the asymptotic expansion comes from the behavior of 
the zeros of $G_2(z,z)$ in the case of $n$-dependent Gaussian
 processes. If $G_2(z,z)$ has a zero of multiplicity $2k$
 on $\partial \D$, i.e., so does the spectral density, 
 then the term of order $(1-r^2)^{-(2k-1)/(2k)}$ appears in
 the asymptotics of $\E N_f(r)$ as $r \to 1$. 
Hence the zeros of the spectral density with the most
 multiplicity determines the asymptotics
 of the second order term. 
Therefore, we obtain the following result for general finitely
dependent cases. 

\begin{cor}\label{cor:general} 
Let $\Xi = \{\xi_k \}_{k \in \Z}$ be the stationary,
 centered, finitely dependent, complex Gaussian process. 
When the spectral density of $\Xi$ has zeros
 $\theta_j$ of multiplicity $2k_j$ for $j=1,2,\dots,p$, 
we set $\alpha = (2k-1)/(2k)$ with $k = \max_{1\le j \le p} k_j$; $\alpha=0$ otherwise. 
Then, there exists a positive constant $C_{\Xi}$ such that 
the expected number of zeros of the GAF $f$ with
 coefficients $\Xi$ within $\D(r)$ is given by 
\[
 \E N_f(r) = \frac{r^2}{1-r^2} - C_{\Xi}
 (1-r^2)^{-\alpha} +o((1-r^2)^{-\alpha}), 
\quad r \to 1. 
\]  
\end{cor}

For example, the Gaussian process $\Xi$ with 
$G_2(z,z) = (const.) \prod_{j=1}^p |z + a_j|^{2k_j}$ 
for $z, a_1,\dots, a_p \in \partial \D$
and $k_1, \dots, k_p \ge 1$  
gives an example of the GAF described in
Corollary~\ref{cor:general}.

This paper is organized as follows. In Section 2, we
recall the Edelman-Kostlan formula and derive its variants
for later use, and prove Theorem~\ref{thm:negative}. 
We also give some examples to give our idea for computation
of the expected number of zeros. 
In Section 3, we prove Theorem \ref{thm:mainintro}. In
Section 4, we briefly recall the method of Puiseux expansion
and prove Theorem \ref{thm:mainintro2}. 

\section{The expected number of zeros: examples} 
\subsection{Expected number of zeros} 
To prove Theorem \ref{thm:mainintro} and
\ref{thm:mainintro2}, we recall the Edelman-Kostlan formula
for the expected number of zeros of GAF. 
\begin{prop}\label{prop:Edelman-Kostlan}
Let $D \subset \C$ be a domain with smooth boundaries, 
$f$ a GAF defined in a neighborhood of $D$, and 
$N_f(D)$ be the number of zeros of $f$ inside $D$. Then,
\begin{align*}
\E N_f(D)
= \frac{1}{4\pi}\int_{D}\Delta\log K_f(z,z)dm(z)
=\frac{1}{2\pi
 \ii}\oint_{\partial D}\partial_z 
\log K_f(z,z)dz,  
\end{align*}
assuming that no singularity lies on
$\partial D$ for the second equality,  
where $dm(z)$ is the Lebesgue measure on the complex plane $\C$ 
and $\ii=\sqrt{-1}$ is the imaginary unit. 
\end{prop}

For the proof of the first equality, see \cite{HKYV}. For
the second equality, the Stokes theorem is used as in \cite{JB,KN}.

In our setting, we have much simpler expressions for $\E N_f(r)$.
\begin{cor}\label{cor:expected_number}
Let $f$ be a GAF defined in \eqref{eq:GAF} with
 \eqref{eq:general_covariance} and \eqref{eq:Gzw}. 
Let $D \subset \D$ be a domain with smooth boundaries and 
$N_f(D)$ be the number of zeros inside $D$. Then, 
\begin{equation}
\E N_f(D)=\frac{1}{2\pi\ii} \oint_{\partial D} \frac{\zbar}{1-|z|^2} dz +\cJ(D), 
\label{eq:ENFD}
\end{equation}
where $\cJ(D)$ has two expressions as follows: 
\begin{equation}
 \cJ(D) 
= \frac{1}{2\pi \ii} \oint_{\partial D}
\frac{G'(z)}{G_2(z,z)}dz
\label{eq:Jr} 
\end{equation}
and 
\begin{equation}
\cJ(D) = -\frac{1}{\pi}\int_{D} 
\left(\frac{|G'(z)|}{G_2(z,z)}\right)^2
 dm(z).
\label{eq:negative}
\end{equation}
In particular, when $D = \D(r)$, \eqref{eq:ENFD} becomes 
\begin{equation}
\E N_f(r)=\frac{r^2}{1-r^2}+\cJ(r),  
\end{equation}
where we simply write $\cJ(r)$ for $\cJ(\D(r))$. 
\end{cor}

\begin{proof}
The first expression \eqref{eq:Jr} directly follows from 
\eqref{eq:general_covariance}, \eqref{eq:Gzw} and the second
 equality in Proposition~\ref{prop:Edelman-Kostlan}. 
For the second expression \eqref{eq:negative}, 
since $\overline{\partial_z G(z)} = \partial_{\zbar} (\overline{G(z)})$, it is easy to see from the first equality in 
Proposition~\ref{prop:Edelman-Kostlan} that
\[
\cJ(D) 
= \frac{1}{\pi} \int_{D}
\partial_z \partial_{\zbar}  \log G_2(z,z) dm(z)
= - \frac{1}{\pi} \int_{D} 
 \frac{|\partial_zG(z)|^2}{(1+G(z)+\overline{G(z)})^2} dm(z). 
\]
This completes the proof. 
\end{proof}

The expression \eqref{eq:negative} essentially, but
 not explicitly, appeared in \cite{SMJM}. 
They derived a similar expression from one-point
 correlation and used to evaluate the expected number of zeros in
 the case of fractional Gaussian noise. 

\begin{rem}
In our setting, $G(z)$ is a polynomial. 
By the change of variables $z \mapsto rz$ in \eqref{eq:Jr} with $D = \D(r)$,
 we have 
\begin{equation}
 \cJ(r) 
= \frac{r}{2\pi \ii} \oint_{\partial \D} 
\frac{G'(rz)}{\Theta(r,z)} dz,  
\label{eq:JR} 
\end{equation}
where $\Theta(r,z)$ is the rational function of $z$ obtained from 
$G_2(rz,rz)$ by putting $\zbar = z^{-1}$ on $\partial \D$. 
In particular, when $\gamma(k)$ is real for every
$k\in\Z$, we have 
\[
 \Theta(r,z) 
= \sum_{k \in \Z} \gamma(k) r^{|k|} z^k. 
\]
Note that $\Theta(1,e^{i\theta})$ is the spectral density at
 least for finitely dependent Gaussian processes. 
Then, one can apply the residue theorem, and from this point
 of view, the behavior of zeros of $\Theta(r,z)$ as $r \to 1$ is
 essential for the order of $\cJ(r)$. 
\end{rem}
Theorem~\ref{thm:negative} is a direct consequence of 
the second expression \eqref{eq:negative} of $\cJ(D)$. 
\begin{proof}[Proof of Theorem~\ref{thm:negative}]
The error term $\cJ(D)$ is clearly non-positive from
 \eqref{eq:negative}. 
Moreover, the right-hand side of \eqref{eq:negative} 
is zero if and only if $G'(z) = 0$ $m$-a.e. $D$. It follows from the uniqueness theorem that $G'(z)$ is identically zero on $\D$, and thus so is $G(z)$ since $G(0)=0$. Therefore, $f(z)$ is equal to $\fpv(z)$ in law. 
\end{proof}

\subsection{Examples} 
In this subsection, we show two examples to see how the
expected number of zeros behaves as $r \to 1$. Although all computations are rather straightforward, they are helpful for understanding of the situation. 
\begin{ex}[Ornstein-Uhlenbeck process]
Let $\gamma(k) = \rho^{|k|} \ (0 < \rho < 1)$. 
The corresponding stationary Gaussian process is the
 (discrete time) Ornstein-Uhlenbeck process. In this case,
 we see that 
$G(z) = \rho z (1-\rho z)^{-1}$ and 
\[
G_2(z,w) 
= \frac{1-\rho^2 z \wbar}{(1-\rho z)(1- \rho
 \wbar)}. 
\]
By using $\zbar = z^{-1}$ for $z \in \partial \D$, we see
 that 
\[
\Theta(r,z) 
= \frac{z(1-\rho^2 r^2)}{(1-\rho r z)(z- \rho r)}. 
\]
We apply \eqref{eq:JR} to this case. The only zero $z=0$ of 
$\Theta(r,z)$, which does not move in $r$, contributes to the
 residue as the only pole. 
Hence, we have 
\[
 \E N_f(r) = \frac{r^2}{1-r^2} - \frac{\rho^2 r^2}{1-\rho^2
 r^2} = \frac{r^2}{1-r^2} - \frac{\rho^2}{1-\rho^2}
+ O(1-r^2), \quad r \to 1.
\]
In this case, $G(z)$ is analytic in $\D(1/\rho)$ and
$\Theta(1,z)$, or equivalently $G_2(z,z)$, does not vanish on $\partial \D$. 
\end{ex}

\begin{rem}
As was seen in this example, the second term $\cJ(r)$ is $O(1)$ as $r \to 1$
 whenever $G(z)$ is analytic in a neighborhood of 
$\overline{\D} := \D \cup \partial \D$ 
and $\Theta(r,z)$ does not vanish on $\partial \D$. 
\end{rem}

\begin{ex} 
For $0<\rho<1$, let $\zeta$ and $\{\eta_k\}_{k \in \Z}$ 
be i.i.d. complex standard normal random variables and
 define the Gaussian process $\Xi = \{\xi_k\}_{k \in \Z}$ by 
\[
\xi_k = \sqrt{\rho} \zeta + \sqrt{1-\rho} \eta_k \quad
 \text{for $k \in \Z$}.  
\]
Then, the corresponding GAF is equal in law to 
\begin{equation}
\sqrt{\rho} \frac{\zeta}{1-z} + \sqrt{1-\rho} \fpv(z)
\label{eq:inlaw} 
\end{equation}
and its covariance function is given by 
\[
 \gamma(k) = 
\begin{cases}
 1 & k=0, \\
 \rho & \text{otherwise}. 
\end{cases}
\]
In this case, $G(z) = \rho z(1-z)^{-1}$ and 
\begin{align*}
G_2(z,z)=\frac{1-(1-\rho)(z+\zbar)+(1-2\rho)|z|^2}{(1-z)(1-\zbar)},
\end{align*}
and hence 
\begin{align*}
\Theta(r,z)
=- \frac{(1-\rho)rz^2 - (1+(1-2 \rho) r^2)z +(1-\rho) r}{(1-rz)(z-r)} 
\end{align*}
The zeros of $\Theta(r,z)$ are $\nu$ and $\nu^{-1}$, where $\nu=\frac{\delta-\sqrt{\delta^2-4}}{2}$ and $\delta=\frac{1+(1-2\rho)r^2}{(1-\rho)r}$. 
Note that $\nu$ (resp., $\nu^{-1}$) is inside (resp., outside) $\D$. 
By using \eqref{eq:JR} and the residue theorem, we have
\begin{align*}
\E
 N_f(r)=\frac{r^2}{1-r^2}-\frac{\rho}{1-\rho}\frac{\nu-r}{(\nu-\nu^{-1})(1-\nu
 r)}. 
\end{align*}
As $r \to 1$, we have 
\[
 \E N_f(r)
= \frac{r^2}{1-r^2} - \frac{1}{2}\sqrt{\frac{\rho}{1-\rho}}
\frac{1}{\sqrt{1-r^2}} + O(1). 
\]
\end{ex}
\begin{rem}
(i) The convergence radius of $G(z)$ is $1$ and its singularity
 is located only at $z=1$. 
The zeros of $\Theta(r,z)$ are $\nu$ and $\nu^{-1}$ given above, where 
$\nu$ (resp., $\nu^{-1}$) is inside (resp., outside) $\D(r)$. 
Both $\nu$ and $\nu^{-1}$ converge to $1$ as $r \to 1$ and 
the second term of $O((1-r^2)^{-1/2})$ comes from
 $(\nu-\nu^{-1})^{-1}$ as the residue at $z=\nu$. \\
(ii) From \eqref{eq:inlaw} we intuitively observe that near
 $z=1$, the first term $\zeta/(1-z)$ pushes up
 the absolute values of $\sqrt{1-\rho} \fpv(z)$ and decreases the
 number of zeros. 
\end{rem}

We would like to emphasize that the behavior of zeros of
$\Theta(r,z)$ 
as $r \to 1$ is essential for the asymptotic behavior of the error term $\cJ(r)$. 

\section{$2$-dependent cases} 
In this section, we prove Theorem \ref{thm:mainintro}. 
\subsection{Case (I)}
First we show Case (I). 
\begin{proof}[Proof of Case (I) in Theorem \ref{thm:mainintro}] 
First we note that $G(z)=az+bz^2$ and then 
$\Theta(r,z) = 1 + ar(z+z^{-1}) + br^2(z^2+z^{-2})$. 
From \eqref{eq:JR}, we have 
\begin{equation}
\mathcal{J}(r)=\frac{r}{2\pi\ii}\oint_{\partial\mathbb{D}}\frac{a+2brz}{1+ar(z+z^{-1})+br^2(z^2+z^{-2})}dz,
\label{eq:Jrcase1}
\end{equation}
We suppose $(a,b) \in \partial \cP_1 \cap \partial \cP$, i.e.,
$a = \pm 2 \sqrt{b(1-2b)}$ with $1/6 \le b \le 1/2$. 
By the symmetry, it is enough to consider the case $a>0$. 
Since the denominator is
 reciprocal, if $\gamma$ is one of its roots, then the roots
 are given as $\gamma, \gamma^{-1}, \gammabar, \gammabar^{-1}$. Here
 we suppose $\gamma \in \D$ and in the upper-half
 plane. Thus, $\gamma, \gammabar$ (resp., $\gamma^{-1},
 \gammabar^{-1}$) are inside (resp., outside) $\D$. 
By taking the residues at $\gamma$ and $\gammabar$, we see
 that 
\begin{align*}
\mathcal{J}(r)
&= \frac{1}{2\pi\ii br}
 \oint_{\partial\D}\frac{z^2(a+2brz)}{(z-\gamma)(z-\gammabar)(z-\gamma^{-1})(z-\gammabar^{-1})} dz
 \\
&= 
\frac{2}{br} \Re \left(
\frac{\gamma^2(a+2br\gamma)}{(\gamma-\gammabar)(\gamma-\gamma^{-1})(\gamma-\gammabar^{-1})}
\right). 
\end{align*}
Let $X=z+z^{-1}$ and rewrite the denominator as
 $br^2 X^2+arX+1-2br^2$, whose roots are distinct and given by 
$X_{\pm}=(-a \pm \ii 2 \sqrt{2}b\sqrt{1-r^2})/(2br)$. 
It is easy to see that 
\begin{align*}
 \gamma = \frac{X_- + \sqrt{X_-^2-4}}{2}, 
\quad \gammabar = \frac{X_+ + \sqrt{X_+^2-4}}{2}, 
\ \gamma^{-1} = \frac{X_- - \sqrt{X_-^2-4}}{2}, 
\quad \gammabar^{-1} = \frac{X_+ - \sqrt{X_+^2-4}}{2}.
\end{align*}
Here we take the branch of $\sqrt{z}$ such that
 $\sqrt{1}=1$ and analytic in $\C \setminus (-\infty,0]$. 
Note that 
\begin{align}
\gamma-\gamma^{-1} = \sqrt{X_{-}^2-4} 
&=\frac{1}{br}\left(\sqrt{\frac{\alpha+\sqrt{\alpha^2+\beta^2}}{2}}+\ii\sqrt{\frac{-\alpha+\sqrt{\alpha^2+\beta^2}}{2}}\right)
\label{eq:imgamma} 
\end{align}
with $\alpha=b-2b^2(r^2+2)$ and
 $\beta=2b\sqrt{2(1-r^2)b(1-2b)}$. 
It is easy to see that
\begin{align*}
(\gamma-\gammabar)(\gamma-\gammabar^{-1})=\gamma(X_--X_+)=-\gamma\frac{2\sqrt{2(1-r^2)}}{r}\ii
\end{align*}
and hence 
\begin{equation}
\cJ(r) 
=-\frac{1}{b\sqrt{2(1-r^2)}}\Im\left(\frac{\gamma(a+2br\gamma)}{\gamma-\gamma^{-1}}\right). 
\label{eq:JrcaseI}
\end{equation}
We note that $\gamma=(X_{-} + \gamma - \gamma^{-1})/2$. 
Substituting it to the numerator and expanding it by $Y := \gamma-\gamma^{-1}$, we have
\begin{align}
\frac{\gamma(a+2br\gamma)}{\gamma-\gamma^{-1}}
&=\frac{1}{2Y} \Big( X_- (a+brX_-) 
 + (a+2brX_-) Y + br Y^2\Big) \nonumber \\ 
&=\frac{2br^2-1}{2r} Y^{-1} - \ii b\sqrt{2(1-r^2)}+\frac{brY}{2}. 
\label{eq:Yexpansion}
\end{align}
Here we used the fact that $X_-$ is a solution of the equation $br^2 X^2 + arX + 1-2br^2=0$. 
Since $\alpha = -b(6b-1) + O(1-r^2)$ and $\beta = 2b\sqrt{2b(1-2b)}\sqrt{1-r^2}$, we see that 
\begin{equation}
 \Im Y = \sqrt{\frac{6b-1}{b}} + O(1-r^2), \quad 
 \Im Y^{-1} = - \sqrt{\frac{b}{6b-1}} + O(1-r^2),\quad r\to1. 
\label{eq:Yasymptotics1}
\end{equation}
Hence it follows from \eqref{eq:JrcaseI}, \eqref{eq:Yexpansion} and \eqref{eq:Yasymptotics1} that 
\begin{align*}
\cJ(r)=-\sqrt{\frac{2b}{6b-1}}\frac{1}{\sqrt{1-r^2}}+O(1),\quad r\to1.
\end{align*}
This completes the proof of Case (I).
\end{proof}

\subsection{Case (II)}
Next we prove Case (II). 
\begin{proof}[Proof of Case (II) in Theorem \ref{thm:mainintro}]
By the symmetry, it is enough to consider the case $b=a-1/2$
 $(-1/2\leq b\leq 1/6)$. We divide the proof of Case (II)
 into two cases, i.e., (i) $0<b\leq 1/6$ and (ii) $-1/2\leq b\leq 0$. 
In this subsection, we always consider the situation for $r$ sufficiently close to $1$
 depending on $b$. 

First we prove the case (i). 
The roots of $br^2X^2+arX+1-2br^2=0$ are real and given by
 $X_{\pm}=(-a \pm \la)/2br\in\R$ 
with $\la = \sqrt{a^2-4b^2 + 8b^2r^2}$. 
Note that $X_\pm^2-4\geq 0$, and 
$X_{+} \to -2$ and $X_{-} \to (2b-1)/(2b)$ as $r \to 1$
As in Case (I), by \eqref{eq:Jrcase1}, since the denominator is
 reciprocal, if two real roots $\gamma$ and $\kappa$ lie inside
 $\D$ such that $\gamma < \kappa < 0$, then all the roots are given as
 $\gamma,\gamma^{-1},\kappa,\kappa^{-1}$. Here
 $\gamma,\kappa$ (resp. $\gamma^{-1},\kappa^{-1}$) are in
 $\D\cap\R$ (resp. in $\D^c\cap\R$), which are given by 
\begin{align}
\gamma=\frac{X_++\sqrt{X_+^2-4}}{2}, \ 
\gamma^{-1}=\frac{X_+-\sqrt{X_+^2-4}}{2}, \
\kappa=\frac{X_-+\sqrt{X_-^2-4}}{2}, \ 
\kappa^{-1}=\frac{X_- - \sqrt{X_-^2-4}}{2}.  
\label{eq:gammakappa}
\end{align}
By \eqref{eq:Jrcase1} and the residue theorem, we see that
\begin{align}
\cJ(r) &= 
\frac{1}{2\pi\ii br }\oint_{\partial\D}\frac{z^2(a+2brz)}{(z-\gamma)(z-\gamma^{-1})(z-\kappa)(z-\kappa^{-1})}dz\nonumber \\
&=\frac{1}{br}\left\{\frac{\gamma^2(a+2br\gamma)}{(\gamma-\gamma^{-1})(\gamma-\kappa)(\gamma-\kappa^{-1})}+\frac{\kappa^2(a+2br\kappa)}{(\kappa-\gamma)(\kappa-\gamma^{-1})(\kappa-\kappa^{-1})}\right\} \nonumber\\
&=
\frac{1}{\la}
\left\{\frac{\gamma(a+2br\gamma)}{\gamma-\gamma^{-1}} - 
\frac{\kappa(a+2br\kappa)}{\kappa-\kappa^{-1}}\right\}.  
\label{eq:Jr21}
\end{align}
Here we used 
\begin{align*}
(\gamma-\kappa)(\gamma-\kappa^{-1})=\gamma(X_+-X_-) = 
\frac{\gamma \la}{br},\quad
 (\kappa-\gamma)(\kappa-\gamma^{-1})=\kappa(X_--X_+) = 
-\frac{\kappa \la}{br}. 
\end{align*}
Since $(\kappa-\kappa^{-1})^{-1}=O(1)$, it suffices to focus on the first term of 
\eqref{eq:Jr21}. We again use the expansion in \eqref{eq:Yexpansion} and have
\begin{align*}
Y = \gamma-\gamma^{-1}
&= 2\sqrt{\frac{1-2b}{1-6b}}\sqrt{1-r^2} + O(1-r^2),\quad r\to1.
\end{align*}
Therefore, 
\begin{align*}
\cJ(r)=-\frac{1}{2}\sqrt{\frac{1-2b}{1-6b}}\frac{1}{\sqrt{1-r^2}}+O(1),\quad r\to1.
\end{align*}
Next we prove the case (ii) of (II). Computation is
 almost the same as in the case (i) of (II), but we only need to change the roles of $\gamma, \gamma^{-1},
 \kappa, \kappa^{-1}$.  
Indeed, $\gamma$ and $\kappa^{-1}$
 (resp. $\gamma^{-1},\kappa$) in \eqref{eq:gammakappa} are
 in $\D\cap\R$ (resp. in $\D^c\cap\R$). By
 \eqref{eq:Jrcase1}, \eqref{eq:gammakappa} and 
\[
 (\kappa^{-1} - \gamma) (\kappa^{-1} - \gamma^{-1}) = \kappa^{-1}(X_--X_+)= -
 \frac{\kappa^{-1} \la}{2br},  
\]
we see that 
\begin{align*}
\cJ(r) 
&=\frac{1}{br}\left\{\frac{\gamma^2(a+2br\gamma)}{(\gamma-\gamma^{-1})(\gamma-\kappa)(\gamma-\kappa^{-1})}+\frac{\kappa^{-2}(a+2br\kappa^{-1})}{(\kappa^{-1}-\gamma)(\kappa^{-1}-\gamma^{-1})(\kappa^{-1}-\kappa)}\right\}\\
&=\frac{2}{\la} 
\left\{\frac{\gamma(a+2br\gamma)}{\gamma-\gamma^{-1}} -  
\frac{\kappa^{-1}(a+2br\kappa^{-1})}{\kappa^{-1}-\kappa}\right\}
\\ 
&=-\frac{1}{2}\sqrt{\frac{1-2b}{1-6b}}\frac{1}{\sqrt{1-r^2}}+O(1),\quad r \to 1. 
\end{align*}
This completes the proof of Case (II).
\end{proof}

\begin{rem}
By the continuity, we have the same asymptotic in Case (II), 
but the behavior of roots $\gamma,\gamma^{-1},\kappa,\kappa^{-1}$ in (II) is
 completely different from Case (I). Indeed,
 $\gamma,\gamma^{-1}\to-1$ and
 $\kappa,\kappa^{-1}\to(2b-1)/4b\pm\sqrt{(1-6b)(1+2b)}/2|b|$
 as $r\to 1$ in Case (II). That is, there is only one pair
 of roots toward the boundary $\partial\D$ as $r\to1$
 except $b=-1/2$. This implies that the asymptotic order is
 affected by the degeneracy of roots of $\Theta(1,z)$ located on the boundary $\partial \D$.
\end{rem}
\subsection{Case (III)}
We give a proof of Case (III). 
\begin{proof}[Proof of Case (III) in Theorem \ref{thm:mainintro}]
Suppose $(a,b) = (2/3, 1/6)$. 
Since $\alpha=\frac{1}{18}(1-r^2)$ and
 $\beta=\frac{1}{9}\sqrt{2(1-r^2)}$, 
by \eqref{eq:imgamma}, we have  
 \begin{align*}
Y =  \gamma-\gamma^{-1} =
  \frac{1}{r}\left(
\sqrt{(1-r^2)+\sqrt{(1-r^2)(9-r^2)}} + 
\ii\sqrt{-(1-r^2)+\sqrt{(1-r^2)(9-r^2)}}
\right).  
 \end{align*}
It easily follows from this expression that 
$\Im Y = O\big( (1-r^2)^{1/4} \big)$ and 
\begin{equation*}
\Im Y^{-1} = - 2^{-7/4} (1-r^2)^{-1/4} + O\big( (1-r^2)^{1/4} \big),\quad r \to 1. 
\end{equation*}
Hence, from \eqref{eq:JrcaseI} and \eqref{eq:Yexpansion}, we can conclude that 
\begin{align*}
\mathcal{J}(r) 
&= -2^{-5/4} (1-r^2)^{-3/4} + O\Big((1-r^2)^{-1/4}\Big), \quad
 r \to 1. 
\end{align*}
This completes the proof of Case (III). 
\end{proof}
\subsection{Case (IV)}
Finally, we give a sketch of the proof of Case (IV). 
Since all zeros of  $\Theta(r,z)$ stay away from $\partial \D$ 
as $r \to 1$ when $(a,b)$ is in the interior of $\cP$, any singularity
contributing to the asymptotic behavior do not appear on the
boundary $\partial\D$, and hence it suffices to consider as
$r$ equals to $1$. Here we only consider the interior of
$\mathcal{P}_1$ and $a>0$. 
We use the same notations in the proof of Case (I).  
In this case, $X_{\pm}=(-a\pm\ii\la(a,b))/(2b)$ with
$\la(a,b)=\sqrt{4b-8b^2-a^2}$ and we see that 
$(\gamma-\overline{\gamma})(\gamma-\overline{\gamma}^{-1}) =
-\gamma b^{-1} \la(a,b) \ii$. Hence, 
\begin{align*}
C(a,b)=-\cJ(1)=\frac{2}{\la(a,b)}\Im\left(\frac{\gamma(a+2b\gamma)}{\gamma-\gamma^{-1}}\right).
\end{align*}
A little more computation shows that  
\begin{align*}
C(a,b) = \frac{\mu(a,b)-(2b-1)}{2\la(a,b)
 \mu(a,b)}\sqrt{4b^2+2b-a^2+2b\mu(a,b)}-1, 
\end{align*}
where $\mu(a,b)=\sqrt{(1+2b)^2-4a^2}$ and that $C(a,b) > 0$
unless $(a,b)=(0,0)$. 
We omit the other cases since we obtain the results just by repeating the similar computation.

\section{Degenerated cases}
In this section, we give a proof of
Theorem~\ref{thm:mainintro2}. 
From \eqref{eq:JR}, we have 
\[
\cJ(r) 
 = \frac{r}{2\pi \ii} \oint_{\partial \D} 
 \frac{G'(z)}{\Theta(r,z)} dz 
= \frac{r}{2\pi \ii} \oint_{\partial \D} 
\frac{p_n(r,z)}{q_n(r,z)} dz 
\]
where $p_n(r,z) = z^n {2n \choose n} G'(w)|_{w=rz}$ and 
\[
 q_n(r,z) := z^n {2n \choose n} \Theta(r,z)
= z^n \sum_{k=-n}^n {2n \choose n+k} r^{|k|} z^k. 
\]
We note from \eqref{eq:G2degenerated} that 
\[
 q_n(1,z) = (z+1)^{2n}. 
\]
To see the asymptotic behavior of $\E N_f(r)$ as $r \to 1$, we
need that of $z(r)$ for $q_n(r,z(r))=0$. 

\subsection{Behavior of the root $z(r)$ as $r \to 1$} 
We first note that $q_n(1,-1)=0$ and $\partial_z q_n(r,z)|_{(r,z)=(1,-1)}=0$. 
Hence, we cannot apply the implicit function theorem in the variable
$z$ to $q_n(r,z)$. Alternatively, we follow a strategy of
using Puiseux series expansion and Newton polygon method
(cf. \cite{WAL}).  

First we note that 
\begin{align*}
\partial_r q_n(r,z) |_{(r,z)=(1,-1)}
&=
2\sum_{k=1}^n k(-1)^{n+k}\binom{2n}{n+k} \\
&=(-1)^{n+1}\frac{n+1}{2n-1}\binom{2n}{n+1} \neq
0. 
\end{align*}
By shifting $(r,z)\to(1-r,z+1)$ in $q_n(r,z)$, we consider
\begin{equation}
Q_n(x,y):=\sum_{l=0}^{2n}\binom{2n}{l}(1-x)^{|l-n|}(y-1)^l.
\label{eq:Qnxy}
\end{equation}
Note that
$Q_n(0,y)=y^{2n}$. Following \cite{WAL}, we denote by $\C\{x,y\}$
(resp., $\C\{x\}$) the ring of convergent power series
defined by two variables $x,y$ (resp., one variable $x$). 
If $f\in\C\{x,y\}$ satisfies $f(0,y)=y^mA(y)$ with $A(0)\neq
0$, then we say $f$ is {\it regular in $y$ of order
$m$} \cite[p.20]{WAL}. In our setting, $Q_n(x,y)$ is
regular in $y$ of order $2n$. 
We can use the following theorem from \cite[p.20, Theorem
2.2.6]{WAL} to guarantee the existence of $2n$ distinct
solutions to the equation $Q_n(x,y)=0$ around $(x,y)=(0,0)$. 
\begin{thm}[\cite{WAL}]
(i) Any equation $f(x,y)=0$ where $f\in\C\{x,y\}$ with $f(0,0)=0$, $f(0,y)\not\equiv 0 $ admits at least one solution of the form $y=g(x^{1/m_1})\in\C\{x\}$.\par
(ii) If $f$ is regular in $y$ of order $m$, and we write $f=UF$ with $U$ a unit and $F$ a monic polynomial of degree $m$ in $y$, there are $m$ such solutions $g_j(x^{1/m_{j}})$, all distinct unless the discriminant of $F$ vanishes identically, and $F(y)\equiv\prod_{j=1}^m\left(y-g_j(x^{1/m_{j}})\right)$.
\end{thm}

For our purpose, we need more explicit form of $g_j$'s so that we
directly perform the Newton polygon method below. 

\begin{figure}[h]
\centering
\begin{tikzpicture}
   
    \draw[style=help lines,step=.5cm] (0,0) grid (4.2,4.2);
    
    \draw[->,thick] (-0.1,0) -- (4.5,0 ) node[anchor=west]{};
    \draw[->,thick] (0,-0.1) -- (0,4.5) node[anchor=south]{};
    
    \foreach \x/\xtext in  {0/0,0.5/1,1/2, 1.5/3, 2/4, 2.5/5, 3/6, 3.5/7, 4/8} 
    \draw[shift={(\x,0)}] (0pt,2pt) -- (0pt,-2pt) node[below] {$\xtext$};
    \foreach \y/\ytext in  {0/0,0.5/1,1/2, 1.5/3, 2/4, 2.5/5, 3/6, 3.5/7, 4/8} 
    \draw[shift={(0,\y)}] (0pt,2pt) -- (0pt,-2pt) node[left] {$\ytext$};

    \draw[thick] (0,4.0) -- (0.5,0);
    \fill (canvas cs:x=0cm,y=4cm) circle (2pt); 
    
    \fill (canvas cs:x=.5cm,y=0cm) circle (2pt);           
    \fill (canvas cs:x=.5cm,y=0.5cm) circle (2pt);
    \fill (canvas cs:x=.5cm,y=1.0cm) circle (2pt);
    \fill (canvas cs:x=.5cm,y=1.5cm) circle (2pt);
    \fill (canvas cs:x=.5cm,y=3.5cm) circle (2pt);
    \fill (canvas cs:x=.5cm,y=4cm) circle (2pt);
    
    \fill (canvas cs:x=1cm,y=0cm) circle (2pt);
    \fill (canvas cs:x=1cm,y=.5cm) circle (2pt);
    \fill (canvas cs:x=1cm,y=1cm) circle (2pt);
    \fill (canvas cs:x=1cm,y=1.5cm) circle (2pt); 
    \fill (canvas cs:x=1cm,y=3cm) circle (2pt);
    \fill (canvas cs:x=1cm,y=3.5cm) circle (2pt);   
    \fill (canvas cs:x=1cm,y=4cm) circle (2pt);
      
    \fill (canvas cs:x=1.5cm,y=0cm) circle (2pt);
    \fill (canvas cs:x=1.5cm,y=.5cm) circle (2pt);
    \fill (canvas cs:x=1.5cm,y=1cm) circle (2pt);
    \fill (canvas cs:x=1.5cm,y=1.5cm) circle (2pt);
    \fill (canvas cs:x=1.5cm,y=2.5cm) circle (2pt); 
    \fill (canvas cs:x=1.5cm,y=3cm) circle (2pt);
    \fill (canvas cs:x=1.5cm,y=3.5cm) circle (2pt);
    \fill (canvas cs:x=1.5cm,y=4cm) circle (2pt);
    
    \fill (canvas cs:x=2cm,y=0cm) circle (2pt);
    \fill (canvas cs:x=2cm,y=.5cm) circle (2pt);
    \fill (canvas cs:x=2cm,y=1cm) circle (2pt);
    \fill (canvas cs:x=2cm,y=1.5cm) circle (2pt);
    \fill (canvas cs:x=2cm,y=2cm) circle (2pt); 
    \fill (canvas cs:x=2cm,y=2.5cm) circle (2pt); 
    \fill (canvas cs:x=2cm,y=3cm) circle (2pt);
    \fill (canvas cs:x=2cm,y=3.5cm) circle (2pt);
    \fill (canvas cs:x=2cm,y=4cm) circle (2pt);
\end{tikzpicture}
\caption{Newton polygon of $Q_n(x,y)$ for $n=4$. A point
 $(r,s)$ is marked when the coefficient $x^r y^s$ of
 $Q_n(x,y)$ is nonzero. }
\label{fig:newton}
\end{figure}
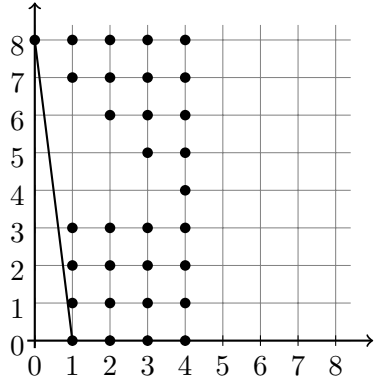

The solution $y(x)$ to $Q_n(x,y)=0$ around the neighborhood
of the origin $(0,0)$ is described by this theorem since
$Q_n(x,y)$ is a bivariate polynomial. 
Now we will compute the asymptotic expansion of $y=y(x)$ in
$Q_n(x,y(x))=0$ at the origin $(0,0)$ following the Newton
polygon method \cite[p.15, Theorem 2.1.1]{WAL}. 
Here we give a brief description of the algorithm following
\cite{WAL}. First, given $f(x,y)=0$, we plot a point $(r,s)$
of exponents for each term $c_{r,s}x^ry^s$ of $f(x,y)$ on
$\R^2$ if $c_{r,s} \not= 0$ and then we have the convex hull
containing all points plotted. Its boundary is made up of
straight line segments which do not lie on the coordinate
axes. It is called the {\it Newton polygon}. Secondly, we
denote by $m_1$ one of the reciprocal numbers of the
negative of a slope among these segments. Then we consider
$f(x,x^{m_1}(a_1+y_1))$ and solve $a_1$ by focusing on the
terms of the lowest degrees in $x$ due to
$f(x,y)=0$. Thirdly, let
$f^{(1)}(x,y_1)=x^{-l} f(x,x^{m_1}(a_1+y_1)$ where $l$ is the
intersection of $s$-axes. Repeat the above process and then
we can obtain the solution
$y=a_1x^{m_1}+a_2x^{m_1+m_2}+\cdots$ of $f(x,y)=0$ for
$f\in\C\{x,y\}$.
For $Q_n(x,y)$, its Newton polygon joins $(1,0)$ and
$(0,2n)$ as shown in Figure~\ref{fig:newton} for $n=4$. 
Thus, it is guaranteed that $Q_n(x,y)=0$ has the solution of
the form 
\begin{align*}
y=x^{1/(2n)}(a_1+y_1), 
\end{align*}
where $y_1=x^{m_2}(a_2+y_2)$ with $m_2 \in \Q$ being positive. 
Setting $t = x^{1/(2n)}$ (equivalently $x=t^{2n}$) in \eqref{eq:Qnxy} for simplicity, we have 
\[
Q(t^{2n}, t(a_1+y_1))
= \sum_{l=0}^{2n}\binom{2n}{l}(1-t^{2n})^{|l-n|}(t(a_1+y_1)-1)^l=0
\]
and the left-hand side can be expanded as follows: 
\begin{align*}
\lefteqn{Q(t^{2n}, t(a_1+y_1))}\\
&=\left(\sum_{l=0}^{2n}\binom{2n}{l}|l-n|(-1)^{l+1}+a_1^{2n}+2na_1^{2n-1}y_1+\binom{2n}{2}
 a_1^{2n-2}y_1^2
 \right)t^{2n} \\
&\quad +\sum_{l=0}^{2n}\binom{2n}{l}|l-n|l(-1)^l(a_1+y_1) t^{2n+1}+\sum_{l=0}^{2n}\binom{2n}{l}|l-n|\binom{l}{2}(-1)^{l-1}(a_1+y_1)^2 t^{2n+2}+O(t^{2n+3})
\end{align*}
Since $y_1=O(x^{m_2}) = O(t^{2n m_2})$ for positive $m \in \Q$, 
the leading term is of order $t^{2n}$ and its coefficient is given by
\begin{align*}
a_1^{2n}+\sum_{l=0}^{2n} \binom{2n}{l}|l-n| 
(-1)^{l+1}
&=a_1^{2n}+2(-1)^{n}\binom{2(n-1)}{n-1}. 
\end{align*}
Thus, $a_1$ is characterized by the solution of the
equation 
\begin{equation}
a_1^{2n}+2(-1)^n\binom{2(n-1)}{n-1}=0.
\label{eq:x2n} 
\end{equation}
For this $a_1$, the term of the lowest order
$t^{2n}$ in $Q_n(t^{2n}, t(a_1+y_1))$ vanishes and we have
\begin{align}
Q_n^{(1)}(t, y_1)
&:=t^{-2n} Q_n(t^{2n},t(a_1+y_1))\label{eq:Qn1ty} \\
&=2na_1^{2n-1}y_1+n(2n-1)a_1^{2n-2}y_1^2 \nonumber \\
&\quad +c(a_1+y_1)t
 +\sum_{l=0}^{2n}\binom{2n}{l}|l-n|\binom{l}{2}(-1)^{l-1}(a_1+y_1)^2 t^2
 + O(t^3), \nonumber
\end{align}
where 
\[
 c 
= \sum_{l=0}^{2n} {2n \choose l} |l-n| l (-1)^l 
= (-1)^{n+1} 2n {2(n-1) \choose n-1} 
\not=0,
\]
which implies $y_1=O(t)$. Now we repeat the same procedure for
$Q_n^{(1)}(t,y_1)$. We substitute $y_1 = t(a_2 + y_2)$ in 
$Q_n^{(1)}(t,y_1)$ and compare the term of order $t$ to
obtain
\[
 ca_1 + 2na_1^{2n-1} a_2= 0, 
\]
and hence 
\begin{equation}
 a_2 
= - \frac{c a_1^{-2(n-1)}}{2n}
= -\frac{1}{2} a_1^2. 
\label{eq:a2}  
\end{equation}

Putting $y_1=t(a_2+y_2)$ in \eqref{eq:Qn1ty} 
and using \eqref{eq:x2n} and \eqref{eq:a2}
yields 
\begin{align*}
t^{-1}Q_n^{(1)}(t, t(a_2+y_2&)) 
= 2na_1^{2n-1}y_2 + \left(c'+cy_2+\cdots\right) t
+O(t^2), 
\end{align*}
where 
\[
c' = n(2n-1)a_1^{2n-2}a_2^2+c a_2 + 
\sum_{l=0}^{2n}\binom{2n}{l}|l-n|{l \choose
2}(-1)^{l-1}a_1^2 \not=0, 
\]
which implies $y_2=O(t)$. 
In summary, by taking \eqref{eq:x2n}, \eqref{eq:a2} and 
$y = t \{a_1 + t (a_2 + O(t))\}$ 
into account, the solutions to the equation $Q_n(x,y)=0$
around $x=0$ are of the form
\begin{equation}
y^{(n)}_j(x) 
= b_j^{(n)} x^{1/(2n)}
-\frac{1}{2} (b_j^{(n)})^2 x^{1/n}  + O(x^{3/(2n)}),
			  \quad \text{as $x \to 0$},
\label{eq:puiseux4y}
\end{equation}
for $j=0,1,\dots,2n-1$, where
$\{b^{(n)}_j\}_{j=0}^{2n-1}$ are the solutions of
\eqref{eq:x2n}.  

\begin{prop}\label{prop:znj}
Let $q_n(r,z)= z^n\sum_{k=-n}^n \binom{2n}{n+k}r^{|k|}z^k$. Then,
 the solutions $z=z^{(n)}_j(r)$ to the equation $q_n(r,z)=0$
 are of the form
\begin{equation}
z^{(n)}_j(r) = -1 + b^{(n)}_j (1-r)^{\frac{1}{2n}} -
\frac{1}{2} (b^{(n)}_j)^2 (1-r)^{\frac{1}{n}} +
 O((1-r)^{\frac{3}{2n}}), \quad r \to 1, 
\label{eq:znj}
\end{equation}
where 
 \begin{equation}
  b^{(n)}_j = \left\{2\binom{2(n-1)}{n-1}\right\}^{1/(2n)} 
 \exp\left(\frac{2 j - n+1}{2n} \pi \ii \right) 
 \quad (j=0,1,\dots,2n-1). 
 \label{eq:bnj}
 \end{equation}
\end{prop} 
\begin{proof}
Since $z_j^{(n)}(r) = -1 + y_j^{(n)}(1-r)$, putting $x=1-r$ and
 $y=z+1$ in \eqref{eq:puiseux4y} yields 
\begin{align*}
 z_j^{(n)}(r) + 1 
&= b_j^{(n)} (1-r)^{\frac{1}{2n}} 
-\frac{1}{2} (b_j^{(n)})^2 (1-r)^{\frac{1}{n}} 
+ O\Big((1-r)^{\frac{3}{2n}}\Big), \\
\end{align*}
as $r \to 1$. 
We obtain the assertion. 
\end{proof}

\subsection{Proof of Theorem~\ref{thm:mainintro2}} 
We first observe the following asymptotics. 
\begin{lem}\label{lem:znj} 
For $k=0,1,\dots,2n-1$, as $r \to 1$, 
\begin{align*}
\prod_{j=0 \atop{j\not=k}}^{2n-1} (z^{(n)}_k(r) - z^{(n)}_j(r)) 
&= 
(2n)(-1)^{n-1} (e_k^{(n)})^{-1} 
\left\{{2(n-1) \choose n-1}\right\}^{\frac{2n-1}{2n}} 
(1-r^2)^{\frac{2n-1}{2n}} \\
&\quad \times \left\{1 - C_n e_k^{(n)} (1-r^2)^{\frac{1}{2n}} + O\Big( (1-r^2)^{\frac{1}{n}} \Big)\right\},  
\end{align*}
where $C_n$ is a constant depending only on $n$ and 
\begin{equation}
e^{(n)}_k 
= 
\exp\left(\frac{2 k - n+1}{2n} \pi \ii \right) 
\quad (k=0,1,\dots,2n-1). 
\label{eq:Cnj0}
\end{equation}
\end{lem}
\begin{proof} From Proposition~\ref{prop:znj}, we have 
\begin{align*}
\prod_{j=0 \atop{j\not=k}}^{2n-1} (z^{(n)}_k(r) - z^{(n)}_j(r)) 
&= 
\prod_{j=0 \atop{j\not=k}}^{2n-1} (b^{(n)}_k - b^{(n)}_j)
 \cdot (1-r)^{\frac{2n-1}{2n}} \nonumber \\
&\quad -\frac{1}{2} \sum_{l=0\atop{l\not=k}}^{2n-1} 
\prod_{j=0 \atop{j\not=k, l}}^{2n-1} (b^{(n)}_k - b^{(n)}_j)
 \cdot \Big\{(b^{(n)}_k)^2 - (b^{(n)}_l)^2)\Big\} \cdot
 (1-r)^{\frac{2n}{2n}} \nonumber \\
&\quad + O\Big( (1-r)^{\frac{2n+1}{2n}} \Big). 
\end{align*}
Since $\displaystyle \prod_{j=0}^{2n-1} (z - e^{\frac{j-k}{n} \pi \ii}) = z^{2n}-1$, by differentiating both sides and putting $z=1$, we obtain $\displaystyle \prod_{j=0 \atop{j\not=k}}^{2n-1} (1 - e^{\frac{j-k}{n} \pi \ii}) = 2n$ 
for every $k =0,1,\dots, 2n-1$. 
Hence, we have 
 \[
 \prod_{j=0 \atop{j\not=k}}^{2n-1} (e_k^{(n)} - e_j^{(n)})
 = (e_k^{(n)})^{2n-1} 
 \prod_{j=0 \atop{j\not=k}}^{2n-1} (1 - e^{\frac{j-k}{n} \pi \ii}) 
= 2n (-1)^{n-1} (e_k^{(n)})^{-1}.  
 \]
and thus, by \eqref{eq:bnj}, 
\[
\prod_{j=0 \atop{j\not=k}}^{2n-1} (b^{(n)}_k - b^{(n)}_j)
= 
\left\{
2 {2(n-1) \choose n-1}
\right\}^{\frac{2n-1}{2n}}
2n (-1)^{n-1} (e_k^{(n)})^{-1}.  
\]
Similarly, 
\begin{align*}
\sum_{l=0\atop{l\not=k}}^{2n-1} 
\prod_{j=0 \atop{j\not=k, l}}^{2n-1} (b^{(n)}_k - b^{(n)}_j)
 \cdot \Big\{(b^{(n)}_k)^2 - (b^{(n)}_l)^2 \Big\} 
&= 2 {2(n-1) \choose n-1} 2n(-1)^{n-1}
(e^{(n)}_{k})^{-1}
\sum_{l=0\atop{l\not=k}}^{2n-1} (e^{(n)}_k + e^{(n)}_l) \\
&=(-1)^{n-1} 8n(n-1) {2(n-1) \choose n-1}. 
\end{align*}
Since $1-r = \frac{1-r^2}{2} + O((1-r^2)^2)$, we obtain the assertion. 
\end{proof}

Now we give a proof of Theorem~\ref{thm:mainintro2}. 
We appeal to \eqref{eq:Jr} to obtain the asymptotic behavior
of $\cJ(r)$. 
First we remark that the constant $b^{(n)}_j$ in 
\eqref{eq:bnj} lies in the
 right-half plane $\{z\in\C : \Re z > 0\}$ for $j=0,1,\dots,n-1$ and  
 the left-half plane $\{z\in\C: \Re z < 0\}$ for $j=n, n+1, \dots, 2n-1$. 
Thus, if $r$ is sufficiently close to $1$, 
$z^{(n)}_j(r)$ for $j=0,1,\dots,n-1$ lie inside $\D$ and $z^{(n)}_j(r)$ for $j=n+1,n+2,\dots,2n-1$ lie outside $\D$. 
Therefore, we have 
\begin{align}
 \cJ(r) 
&= \frac{r}{2\pi \ii} \oint_{\partial \D}
 \frac{p_n(r,z)}{q_n(r,z)}dz \nonumber \\
&= r 
\sum_{k=0}^{n-1} 
\res\left( \frac{p_n(r,z)}{\prod_{j=0}^{2n-1} (z -
 z^{(n)}_j(r))} ; z= z^{(n)}_k(r)\right)\nonumber \\
&= r \sum_{k=0}^{n-1} 
\frac{p_n(r,z^{(n)}_k(r))}{\prod_{j=0, j \not=k}^{2n-1}
 (z^{(n)}_k(r) - z^{(n)}_j(r))}. 
\label{eq:Jrfinal}
\end{align}
Since $p_n(1,-1) = (-1)^n {2(n-1) \choose n-1}$, from Lemma~\ref{lem:znj} and 
\[
p_n(r, z_k^{(n)}(r)) = 
p_n(1,-1) \left\{1 + 
C_n' e_k^{(n)} (1-r^2)^{1/(2n)} + O\big( (1-r^2)^{1/n} \big) 
\right\}, 
\]
we have 
\begin{align*}
\frac{p_n(r,z^{(n)}_k(r))}{\prod_{j=0, j \not=k}^{2n-1}
 (z^{(n)}_k(r) - z^{(n)}_j(r))}
&= 
\frac{-1}{2n}{2(n-1) \choose n-1}^{\frac{1}{2n}} e_k^{(n)}
(1-r^2)^{-\frac{2n-1}{2n}} \\
&\quad \times 
\left\{1 + (C_n+C_n') e_k^{(n)} (1-r^2)^{\frac{1}{2n}} +
 O\big( (1-r^2)^{\frac{2}{2n}} \big)\right\}, 
\end{align*}
where $C_n'$ is a constant depending only on $n$. 
It is easy to see that 
\begin{equation}
 \sum_{k=0}^{n-1} e_k^{(n)} = (\sin\frac{\pi}{2n})^{-1}, \quad  
 \sum_{k=0}^{n-1} (e_k^{(n)})^2 = 0. 
\label{eq:cancel} 
\end{equation}
Therefore, from \eqref{eq:Jrfinal}, we obtain 
\[
\cJ(r) = 
\frac{-1}{2n \sin(\frac{\pi}{2n})}{2(n-1) \choose n-1}^{\frac{1}{2n}} 
(1-r^2)^{-\frac{2n-1}{2n}} \Big( 1+  O\big(
(1-r^2)^{\frac{2}{2n}} \big) \Big). 
\]
This completes the proof. 

\begin{rem}\label{rem:cancel} 
 A naive computation gives only the error term $O\big( (1-r)^{-(n-1)/n}
 \big)$. Here we saw the cancellation as the second
 equality in \eqref{eq:cancel} to obtain $O\big(
 (1-r)^{-(2n-3)/(2n)} \big)$, which matches the direct computation
 in Case (III) for $n=2$. 
\end{rem}
\begin{rem}
This method can be applied to all cases of finitely dependent Gaussian processes. 
Indeed, the zero of $\Theta(1,e^{\ii\theta})$ of order $2k$ contributes to $\cJ(r)$ as 
constant multiple of $(1-r^2)^{-\frac{2k-1}{2k}}$. 
\end{rem}

\section*{Acknowledgments}
This work was supported by JSPS KAKENHI Grant Number (B) JP18H01124.
KN was also supported by the WISE program (JSPS) at Kyushu University.
TS was also supported by JSPS KAKENHI Grant Numbers
JP16H06338, JP20H00119 and JP20K20884.


\end{document}